\documentclass[11pt,twoside]{amsart}

\usepackage{amssymb}
\usepackage{latexsym}

\newcommand{\Id}{\mathrm{Id}}

\newcommand{\dopu}{{:}\allowbreak\ }
\newcommand{\eps}{\varepsilon}
\newcommand{\cal}{\mathcal}

\newcommand{\mytilde}[1]{\mathbin{\tilde#1}}

\newcommand{\loglike}[1]{\mathop{\rm #1}\nolimits}
\newcommand{\ex}{\loglike{ex}}
\newcommand{\stexp}{\loglike{stexp}}

\newcommand{\co}{\loglike{co}}
\newcommand{\coq}{\loglike{\overline{co}}}
\newcommand{\diam}{\loglike{diam}}
\newcommand{\dist}{\loglike{dist}}
\newcommand{\supp}{\loglike{supp}}

\newcommand{\van}{\loglike{van}}

\newcommand{\SDt}{\loglike{{\cal S}\!{\cal D}}}
\newcommand{\narr}{\loglike{{\cal N}\!\!{\cal A}\!{\cal R}}}

\newcommand{\Nar}{\narr}
\newcommand{\SD}{\SDt}

\newcommand{\N}{{\mathbb N}}
\newcommand{\R}{{\mathbb R}}

\newcommand{\calD}{\mathcal{D}}

\theoremstyle{plain}
\newtheorem{thm}{Theorem}[section]

\newtheorem{prop}[thm]{Proposition}
\newtheorem{cor}[thm]{Corollary}
\newtheorem{lemma}[thm]{Lemma}
\newtheorem{sublemma}[thm]{Sublemma}

\theoremstyle{definition}
\newtheorem{definition}[thm]{Definition}

\theoremstyle{remark}

\numberwithin{equation}{section}

\newcommand{\norm}[1]{\left\Vert#1\right\Vert}

\newcommand{\rest}[2]{#1\raisebox{-0.3ex}{\mbox{$\mid_{#2}$}}}
\def\DP{Daugavet property}
\newcommand{\sDo}{strong Daugavet operator}

\newcommand{\begsta}{\begin{statements}}
\newcommand{\begaeq}{\begin{aequivalenz}}
\def\endsta{\end{statements}}
\def\endaeq{\end{aequivalenz}}
\newcommand{\bea}{\begin{eqnarray*}}
\newcommand{\eea}{\end{eqnarray*}}
\newcommand{\kref}[1]{(\ref{#1})}

\newcounter{abc}   
\newcounter{iiiii} 

\newenvironment{aequivalenz}
{\setcounter{iiiii}{0}
\begin{list}%
{{\rm (\roman{iiiii})}}
{\usecounter{iiiii}
\parsep=0pt plus 1pt
\topsep=1pt plus 2pt minus 1pt
\itemsep=1pt plus 2pt minus 1pt
\leftmargin=3\baselineskip
\labelsep=.6\baselineskip
\labelwidth=2.4\baselineskip
\rightmargin 0pt}%
}%
{\end{list}}

\newenvironment{statements}%
{\setcounter{abc}{0}
\begin{list}%
{{\rm (\alph{abc})}}
{\usecounter{abc}
\parsep=0pt plus 1pt
\topsep=1pt plus 2pt minus 1pt
\itemsep=1pt plus 2pt minus 1pt
\leftmargin=3\baselineskip
\labelsep=.6\baselineskip
\labelwidth=2.4\baselineskip
\rightmargin 0pt}%
}%
{\end{list}}

\begin{document}

\title { Narrow operators on vector-valued sup-normed spaces} 

\author[D.~Bilik, V.~Kadets, R.~Shvidkoy, G.~Sirotkin, D.~Werner]
{ Dmitriy Bilik, Vladimir Kadets, Roman Shvidkoy, Gleb Sirotkin
  and Dirk Werner }

\address{Faculty of Mechanics and Mathematics, Kharkov National
University,\linebreak
 pl.~Svobody~4,  61077~Kharkov, Ukraine}

\address{Faculty of Mechanics and Mathematics, Kharkov National
University,\linebreak
 pl.~Svobody~4,  61077~Kharkov, Ukraine}
\email{vishnyakova@ilt.kharkov.ua}

\curraddr{Department of Mathematics, Freie Universit\"at Berlin,
Arnim\-allee~\mbox{2--6}, D-14\,195~Berlin, Germany}
\email{kadets@math.fu-berlin.de}

\address{Department of Mathematics, University of Missouri,
Columbia MO 65211}
\email{shvidkoy\_r@yahoo.com}

\address{Department of Mathematics, Indiana University -- Purdue
  University Indianapolis, 402 Backford
  Street, Indianapolis, IN 46202}
\email{syrotkin@math.iupui.edu}

\address{Department of Mathematics, Freie Universit\"at Berlin,
Arnimallee~2--6, \qquad {}\linebreak D-14\,195~Berlin, Germany}
\email{werner@math.fu-berlin.de}

\thanks{The work of the second-named author
was supported by a grant from the {\it Alexander-von-Humboldt
Stiftung}.}

\subjclass[2000]{Primary 46B20; secondary 46B04, 46B28, 46E40, 47B38}

\keywords{Daugavet property, narrow operator, strong Daugavet
  operator, USD-nonfriendly spaces, $C(K,E)$-spaces}


\begin{abstract}
We characterise narrow and strong Daugavet operators on
$C(K,E)$-spaces; these are in a way the largest sensible classes of
operators for which the norm equation $\|\Id+T\| = 1+\|T\|$ is valid.
For certain separable range spaces $E$ including all finite-dimensional ones
and locally uniformly convex ones we show that an unconditionally
pointwise convergent sum of narrow operators 
on $C(K,E)$ is narrow, which implies
for instance the known result that these spaces do not have
unconditional FDDs. In a different vein, we construct two narrow
operators on $C([0,1],\ell_1)$ whose sum is not narrow. 
\end{abstract}

\maketitle

\section{Introduction and preliminaries}

This paper is a follow-up contribution to our paper \cite{KadSW2}
where we defined and investigated narrow operators on Banach spaces
with the \DP. We shall first review some definitions and results from
\cite{KadSSW} and \cite{KadSW2} before we describe the contents of
the present paper.

A Banach space $X$ is said to have the \textit{\DP}
if every
rank-$1$ operator $T \dopu X\to X$ satisfies
\begin{equation}\label{eq1.1}
\|\Id+T \|= 1+\|T\|.
\end{equation}
For instance, $C(K)$ and $L_{1}(\mu)$ have the \DP\ provided that $K$
is perfect, i.e., has no isolated points, and $\mu$ does not have any atoms.
We shall have occasion to
use the following characterisation of the Daugavet
property from \cite{KadSSW}; the equivalence of (ii) and (iii) results
from the Hahn-Banach theorem.

\begin{lemma} \label{l:MAIN}
The following assertions are equivalent:
\begin{aequivalenz}
\item
     $X$ has the Daugavet property.
\item
     For every $x\in S(X)$, $x^*\in S(X^*)$  and $\eps >0$ there
         exists some $y\in S(X)$  such that $x^*(y)>1-\eps$  and
$\norm{x+y}>2-\eps$.
\item
For all $x\in S(X)$ and $\eps>0$, $B(X) = \coq \{z\in B(X)\dopu
\|x+z\|>2-\eps\}$.
\end{aequivalenz}
\end{lemma}

It is shown in \cite{KadSSW} and \cite{Shv1} that (\ref{eq1.1})
automatically extends to wider classes of operators, e.g., weakly
compact ones and, more generally,  those that do not fix copies 
of $\ell_{1}$ or strong Radon-Nikod\'ym operators. (A 
strong Radon-Nikod\'ym operator maps the unit ball into a set with
the Radon-Nikod\'ym property.)
In \cite{KadSW2} we found new proofs of these results based on the
notions of a \sDo\ and a narrow operator. An operator $T\dopu X\to Z$ 
is said to be a {\em strong
Daugavet operator}\/ if for every two elements $x, y \in S(X)$, the
unit sphere of $X$, and
for every $\eps > 0$ there is an element $u \in B(X) $, the unit ball
of $X$,  such that
$\|x+u\|> 2 - \eps$ and $\|T(y-u)\|<\eps$. 
It is almost obvious that a \sDo\ $T\dopu X\to X$ satisfies
(\ref{eq1.1}), and the nontrivial task is now to find sufficient
conditions on $T$ to be strongly Daugavet. In this vein we could show
that for instance strong Radon-Nikod\'ym operators and operators not
fixing copies of $\ell_{1}$ are indeed \sDo s. 

For some applications the concept of a \sDo\ is somewhat too wide.
Therefore we defined an operator $T\dopu X\to Z$ 
to be \textit{narrow}
 if for every two elements $x, y \in S(X)$, every $x^* \in X^* $ 
and every $\eps> 0$ there is an element $u \in B(X) $
such that
$\|x+u\|> 2 - \eps$ and $\|T(y-u)\| + |x^*(y-u)|<\eps$.
It follows that $X$ has the \DP\  if and only if all rank-$1$
operators are \sDo s if and only if there is at least one narrow
operator on $X$.
We denote the set of
all strong Daugavet operators on $X$ by $\SD(X)$ 
and the set of all narrow operators on $X$ by $\Nar(X)$.
Actually, in \cite{KadSW2} we took a slightly different point of view
in that we declared two operators $T_{1}\dopu X\to Z_{1}$ and
$T_{2}\dopu X\to Z_{2}$ to be equivalent if $\|T_{1}x\|= \|T_{2}x\|$
for all $x\in X$; $\SD(X)$ and $\Nar(X)$ should really denote the
sets of corresponding equivalence classes. However, in this paper we
shall not make this point explicitly.

In this paper we shall continue our investigations of this type of
operator, mostly in the setting of vector-valued function spaces
$C(K,E)$. One of the drawbacks of the definition of a \sDo\ is that
the sum of two such operators need not be a \sDo\ whereas the
definition of a narrow operator has some built-in additivity quality.
It remained open in \cite{KadSW2} whether the sum of any two narrow
operators is always narrow, although we could prove this to be true
on $C(K)$, and in general we showed that the sum of a narrow operator
and an operator not fixing $\ell_{1}$ is narrow and that the sum of a
narrow operator and a strong Radon-Nikod\'ym operator is narrow. 
(Note that the sum of two strong Radon-Nikod\'ym operators need not
be a strong Radon-Nikod\'ym operator \cite{Schach}.) Our work in
Section~\ref{sec3}, where we completely characterise strong Daugavet
and narrow operators on $C(K,E)$, enables us to give counterexamples
to the sum problem. 

For this we employ a special feature of $\ell_{1}$ explained in
Section~\ref{sec2}. This section introduces a class of Banach spaces
called \textit{USD-nonfriendly} spaces that are sort of remote from
spaces with the \DP; USD stands for uniformly strongly Daugavet. All
finite-dimensional and all locally uniformly convex spaces fall
within this category, but we haven't been able to decide whether a
reflexive  space must be USD-nonfriendly.

The class of USD-nonfriendly spaces is custom-made for our
applications in Section~\ref{sec4}  where we study pointwise
unconditionally convergent series $\sum_{n=1}^\infty T_{n}$ of narrow
operators on $C(K,E)$. If  $E$  is separable and
USD-nonfriendly, we prove that the
sum operator must be narrow again, which is new even in the case
$E=\R$. To achieve this we take a detour investigating the related
class of $C$-narrow operators following ideas from \cite{KadPop}. An
obvious corollary is the result from \cite{KadPop}
that the identity on $C(K)$ is not a pointwise
unconditional sum of narrow operators, which implies
that $C(K)$ does not admit an unconditional Schauder decomposition
into spaces not containing $C[0,1]$. 

We finish this introduction with a technical reformulation of the
definition of a \sDo. Let
$$
D(x,y,\eps) = \{z\in X\dopu \|x+y+z\|>2-\eps,\ \|y+z\|<1+\eps \}
$$
and
\bea
\calD(X) &=&
\{ D(x,y,\eps)\dopu x\in S(X),\ y\in S(X),\ \eps>0 \}, \\
\calD_{0}(X) &=&
\{ D(x,y,\eps)\dopu x\in S(X),\ y\in B(X),\ \eps>0 \}.
\eea
It is easy to see that $T\dopu X\to Z$ is a \sDo\ if and only if $T$
is not bounded from below on any $D\in \calD(X)$
\cite[Prop.~3.4]{KadSW2}. In Section~\ref{sec3} it will be more
convenient to work with $\calD_{0}(X)$ instead; therefore we
formulate a lemma saying that this doesn't make any difference.

\begin{lemma} \label{lem1.2}
An operator $T\dopu X\to Z$ is a \sDo\ if and only if $T$
is not bounded from below on any $D\in \calD_{0}(X)$.
\end{lemma}

\begin{proof}
We have to show that $T\in \SD(X)$ is not bounded from below on
$D(x,y,\eps)$  whenever $\|x\|=1$, $\|y\|\le 1$, $\eps>0$. By the
above, $T$ is not bounded from below on $D(x,-x,1)$; hence, given
$\eps'>0$, for some $\zeta\in S(X)$ we have $\|T\zeta\|<\eps'$. Now
pick $\lambda \ge0$ such that $y+\lambda\zeta \in S(X)$; then there
is some $z'\in X$ such that
$$
\|x+(y+\lambda\zeta)+z'\|>2-\eps, \ 
\|(y+\lambda\zeta)+z'\| <1+\eps, \ 
\|Tz'\|<\eps';
$$
i.e., $z:= \lambda\zeta + z' \in D(x,y,\eps)$ and $\|Tz\|<3\eps'$. 
\end{proof}

\section{USD-nonfriendly spaces}
\label{sec2}

In this section we introduce a class of Banach spaces that are
geometrically opposite to spaces with the \DP. These spaces will
arise naturally in Section~\ref{sec4}.

\begin{prop}\label{prop5.1}
The following conditions for a Banach space $E$ are equivalent.
\begin{enumerate}
\item \label{SD1}
$\SDt(E)=\{0\}$. 
\item \label{SD3}
No nonzero linear functional on $E$ is a strong Daugavet operator.
\item\label{SD4}
For every $x^*\in S(E^*)$ there exist some $\delta>0$ and
$D\in
\calD(E)$ such that $|x^*(z)|>\delta$ for all $z\in D$.
\item \label{SD2}
Every closed absolutely convex 
subset $A \subset E$ such that for every $\alpha > 0$
and every $D \in {\cal D}(E)$ the intersection $(\alpha A) \cap D$
is nonempty coincides with the whole space $E$. 
\end{enumerate}
\end{prop}

\begin{proof}
The  implications~\kref{SD1} $\Rightarrow $ \kref{SD3} $\Rightarrow $ 
\kref{SD4}  are evident.

\kref{SD4} $\Rightarrow $ \kref{SD2}:
Assume there is some closed absolutely convex 
subset $A\subset E$ with the property from~\kref{SD2} 
that does not
coincide with the whole space $E$. By the Hahn-Banach theorem there is
a functional $x^* \in S(E^*)$ and a number $r>0$
such that  $|x^*(a)| \le r$ for every 
$a \in A$. If $\delta>0$ and $D\in \calD(E)$ are arbitrary,
pick $z\in (\frac\delta r A) \cap D$; this intersection is
nonempty by assumption on $A$. It follows that $|x^*(z)| \le \delta$,
hence \kref{SD4} fails.

\kref{SD2} $\Rightarrow $ \kref{SD1}:
Suppose $T\in \SDt(E)$ and
put $A = \{e \in E\dopu  \|Te\| \le 1\}$. By the definition of  a strong 
Daugavet operator this $A$ satisfies~\kref{SD2}. So  $A=E$ and hence 
$T=0$.
\end{proof}

This proposition suggests the following definition.

\begin{definition}
A Banach space $E$ is said to be an \textit{SD-nonfriendly space} 
(i.e., strong Daugavet-nonfriendly) if $\SD(E) = \{0\}$. 
A space $E$ is said to be a \textit{USD-nonfriendly space} 
(i.e., uniformly strong 
Daugavet-nonfriendly) if there exists an $\alpha > 0$ such that 
every closed absolutely convex  subset $A \subset E$ which intersects
all the elements of  ${\cal D}(E)$ contains $\alpha B(E)$. The
largest admissible $\alpha$
is called the \textit{USD-parameter} of $E$.
\end{definition}

Proposition~\ref{prop5.1} shows that a
USD-nonfriendly space is indeed SD-non\-friendly; but the converse is
false as will be shown shortly.
Also, SD-nonfriendliness is opposite to the \DP\ in that the latter is
equivalent to the condition that every functional is a strong
Daugavet operator.

To further motivate the uniformity condition in the above definition,
we supply a lemma.

\begin{lemma}\label{lemma5.1a}
A Banach space $E$ is USD-nonfriendly if and only if
\begin{enumerate}
\item[($3^*$)]
There exists some $\delta>0$ 
such that for every $x^*\in S(E^*)$ 
there exists $D\in
\calD(E)$ such that $|x^*(z)| >\delta$ for all $z\in D$.
\end{enumerate}
\end{lemma}

\begin{proof}
It is enough to prove the implications (a) $\Rightarrow$ (b)
$\Rightarrow$ (c) for the following assertions about a fixed number
$\delta>0$:
\begsta
\item
There exists a closed absolutely 
convex set $A\subset E$ not containing $\delta
B(E)$ that intersects all $D\in \calD(E)$.
\item
There exists a functional $x^*\in S(E^*)$ such that for all $D\in
\calD(E)$ there exists $z_{D}\in D$ satisfying $|x^*(z_{D})| \le
\delta$.
\item
There exists a closed absolutely 
convex set $A\subset E$ not containing $\delta'
B(E)$ for any $\delta'>\delta$ that intersects all $D\in \calD(E)$.
\endsta

To see that (a) implies (b), pick $u\notin A$, $\|u\|\le\delta$. By
the Hahn-Banach theorem we can separate $u$ from $A$ by means of a
functional $x^*\in S(E^*)$; then we shall have for some number $r>0$
that $|x^*(z)|\le r$ for all $z\in A$ and $x^*(u) >r$. On the other
hand, $x^*(u)\le \|x^*\|\,\|u\|\le \delta$; hence (b) holds for $x^*$.

If we assume (b), we define $A$ to be the closed absolutely 
convex hull of the
elements $z_{D}$, $D\in \calD(E)$, appearing in~(b). 
Obviously $A$ intersects each $D\in  \calD(E)$.
If $\delta' B(E)
\subset A$ for some $\delta'>0$, then since $|x^*| \le \delta $ on $A$,
we must have $|x^*| \le \delta $ on $\delta' B(E)$, i.e., $\delta'\le
\delta$. Therefore, $A$ works in~(c).
\end{proof}

In Proposition~\ref{prop5.1} and Lemma~\ref{lemma5.1a} we may replace
$\calD(E)$ by $\calD_{0}(E)$.

We now turn to some examples.

\begin{prop}\label{propc0}
\mbox{}
\begsta
\item
The space $c_{0}$ is SD-nonfriendly, but not USD-nonfriendly.
\item
The space $\ell_{1}$ is not SD-nonfriendly and hence not USD-nonfriendly
either.
\endsta
\end{prop}

\begin{proof}
(a) 
Theorem~3.5 of \cite{KadSW2} implies that $Te_{k}=0$ for every unit
basis vector $e_{k}$ if $T\in \SDt(c_{0})$. [Actually, the theorem
quoted is formulated for operators on $C(K)$ for compact $K$, but the
theorem works likewise on $C_{0}(L)$ with $L$ locally compact.] Hence
$T=0$ is the only strong Daugavet operator on $c_{0}$. (Another way
to see this is to apply Corollary~\ref{cor3.2a}.)

\def\direktesargument{
Assume $x^*=(x^*(n))\in S(\ell_{1})$ is a strong Daugavet operator.
Without loss of generality we may suppose that $\gamma:= -x^*(1)>0$;
note that an operator $T$ is a strong Daugavet operator if and only
if $-T$ is. Let $\eps<\gamma/2$ and fix $y\in S(c_{0})$ with $x^*(y)
>1-\eps$. By definition of a strong Daugavet operator there is
some $u\in c_{0}$       such that 
$$
\|e_{1}+u\| >2-\eps, \quad
\|u\|<1+\eps, \quad
|x^*(u) - x^*(y)| < \eps,
$$
hence
$$
u(1) > 1-\eps, \quad x^*(u) > 1-2\eps.
$$
On the other hand,
\bea
x^*(u) &=&
x^*(1) u(1) + \sum_{n=2}^\infty x^*(n) u(n) \\
&\le&
-\gamma (1-\eps) + \sum_{n=2}^\infty |x^*(n)| (1+\eps) \\
&=&
-\gamma (1-\eps) + (1-\gamma) (1+\eps) \\
&=&
1-2\gamma + \eps < 1-3\eps
\eea
by our choice of $\eps$; hence we have reached a contradiction. This
shows that $x^*=0$  is the only strong Daugavet functional.
}   

To show that $c_{0}$ is not USD-nonfriendly we shall exhibit a closed
absolutely 
convex set $A$ intersecting each $D\in \calD(c_{0})$, yet containing
no ball. Let $A=2B(\ell_{1}) \subset c_{0}$, i.e.,
$$
A= \biggl\{ (x(n))\in c_{0}\dopu \sum_{n=1}^\infty |x(n)| \le2
\biggr\},
$$
which is closed in $c_{0}$.
Fix $x\in S(c_{0})$ and $y\in S(c_{0})$. If $|x(k)|=1$, say $x(k)=1$,
pick $|\beta|\le2$ such that $y(k)+\beta = 1$. Then $\beta e_{k}\in
D(x,y,\eps) \cap A$ for every $\eps>0$. Obviously, $A$ does not
contain a multiple of $B(c_{0})$.

(b) We claim that $x_{\sigma}^*(x)= \sum_{n=1}^\infty \sigma_{n}x(n)$
defines a strong Daugavet functional on $\ell_{1}$ whenever $\sigma $
is a sequence of signs, i.e., if $|\sigma_{n}|=1$ for all~$n$. 
Indeed, let $x\in
S(\ell_{1}) $, $y\in S(\ell_{1}) $ and $\eps>0$. Pick $N$ such that
$\sum_{n=1}^N |x(n)| > 1-\eps$ and define $u\in S(\ell_{1})$ by
$u(n)=0$ for $n\le N$ and $u(n)= \sigma_{n-N}y(n-N)/ \sigma_{n}$ 
for $n>N$. Then $x^*(u)=
x^*(y)$ and $\|x+u\| > 2-\eps$; hence $z:=u-y\in D(x,y,\eps)$ and
$x^*(z)=0$.
\end{proof}

Next we wish to give some examples of USD-nonfriendly spaces.
Recall that a point of local uniform rotundity of the unit sphere of
a Banach space $E$ (a LUR-point) is a point $x_{0}\in S(E)$ such that
$x_{n}\to x_{0}$ whenever $\|x_{n}\|\le1$ and $\|x_{n}+x_{0}\|\to 2$.

\begin{prop}
If the unit sphere of $E$ contains a LUR-point, then $E$ is a
USD-nonfriendly space with USD-parameter~${\ge 1}$.
\end{prop}

\begin{proof}
Let $x_0 \in S(E)$ be  a LUR-point and $A \subset E$ be a closed
absolutely convex  
subset which intersects all the elements of  ${\cal D}(E)$. In
particular for every fixed $y \in S(E)$ the set $A$ intersects all the 
sets $D(x_0 ,y, \eps) \subset E$, $\eps > 0$. By definition of  a
LUR-point this means that all the points of the form $x_0 - y$, 
$y \in S(E)$, belong to $A$, i.e., $B(E)+x_0 \subset A$. But  $-x_0$ 
is also a LUR-point, so $B(E)-x_0 \subset A$, and by convexity of $A$, 
$B(E) \subset A$.
\end{proof}

\begin{cor}
Every locally
uniformly convex space is USD-nonfriendly with USD-parameter~$2$. 
In particular, the
spaces $L_{p}(\mu)$ are USD-nonfriendly for $1<p<\infty$.
\end{cor}

\begin{proof}
This follows from the previous proposition; that the USD-parameter is
$2$ is a consequence of $B(E)+ x_{0} \subset A$ for all $x_{0}\in
S(E)$; see the above proof.
\end{proof}

It is clear that no finite-dimensional space enjoys the \DP, but more
is true.

\begin{prop}\label{prop2.7}
Every finite-dimensional Banach space $E$ is a USD-non\-friendly
space. 
\end{prop}

\begin{proof}
Assume to the contrary that there is a   
finite-dimensional space $E$ that is not
USD-nonfriendly. 
By Lemma~\ref{lemma5.1a} we can find a sequence of functionals
$(x_{n}^*)\subset S(E^*)$ such that $\inf_{z\in D} |x_{n}^*(z)| \le 1/n$
for each $D\in \calD(E)$. By  compactness of the ball we can pass to
the limit and obtain a functional $x^*\in S(E^*)$ with the property
that $\inf_{z\in D} |x^*(z)| = 0$
for each $D\in \calD(E)$.

Denote $K=\{e \in B(E)\dopu  x^*(e)=1 \}$; this  is a norm-compact
convex set. Let $x_{0}\in K$ be an arbitrary point. 
If we apply the above property to $D(x_0 , -x_0, \eps)$ for
all $\eps > 0$, we obtain, again by compactness,
some $z_{0}$ 
 such that $\|z_{0}-x_0\|=1$, $\|z_{0}\|=2$ and $x^*(z_{0})=0$. We have
$x^*(x_{0}-z_{0})=1$, 
so $x_0-z_{0} \in K$. Therefore
$$
2 \ge \diam K \ge \sup_{y\in K} \|x_{0}-y\| \ge \|x_{0}- (x_{0}-z_{0})\|
=\|z_{0}\|=2;
$$
hence $\diam K=2$ and $x_{0}$ is a diametral point of $K$, meaning
$$
\sup_{y\in K}   \|x_{0}-y\| = \diam K.
$$
But any compact convex set of positive diameter contains a
nondiametral point \cite[p.~38]{Die-LNM}; thus we have reached a
contradiction.
\end{proof}

We shall later estimate the worst possible USD-parameter 
of an $n$-di\-men\-sional normed space.

We haven't been able to decide whether 
every reflexive space is USD-non\-friendly.
Proposition~\ref{expo2} below presents a necessary condition a
hypothetical reflexive USD-friendly ($=$ not USD-nonfriendly)
space must fulfill.

First an easy geometrical lemma.

\begin{lemma}
Let $x,h \in E$, $\|x\| \le 1 + \eps$, $\|h\| \le 1 + \eps$, 
$\|x+h\| \ge 2 - \eps$. Let $f \in S(E^*)$ be a supporting functional 
of $(x+h)/{\|x+h\|}$. Then $f(x)$ as well as $f(h)$ are estimated 
from below by $1 - 2 \eps$.
\end{lemma}

\begin{proof}
Denote $a=f(x)$, $b=f(h)$. Then $\max(a,b) \le 1 + \eps $ but 
$a+b \ge 2 - \eps$. So  $\min(a,b)= a+b - \max(a,b) \ge 1 - 2 \eps $.
\end{proof}

Let $E$ be a reflexive space, $x_0^*$ be a strongly exposed point
of $S(E^*)$ with strongly exposing evaluation functional $x_{0}$; i.e.,
the diameter of the slice $\{x^* \in S(E^*)\dopu  x^*(x_0) > 1 - \eps\}$ 
tends to $0$ when $\eps$ tends to~$0$. 
Denote 
$$
S_{x_0^*}=\{x \in S(E)\dopu  x_0^*(x) = 1\}.
$$

\begin{prop} \label{expo}
Let  $E$, $x_0^*$,  $x_0$ be as above,
$A$ be a closed convex set which intersects all the 
sets $D(x_0 ,0, \eps)$, $\eps > 0$. Then 
$A$ intersects $S_{x_0^*}$.
\end{prop}

\begin{proof}
For every $n \in \N$ select $h_n \in A \cap D(x_0 ,0, \frac{1}{n})$.
Then $\|h_n\| \le 1 + \frac{1}{n}$, $\|x_0+h_n\| \ge 2-\frac{1}{n}$. 
Denote by   $f_n$  a supporting functional 
of ${(x_0+h_n)}/{\|x_0+h_n\|}$. By the previous lemma 
$f_n(x_0)$ tends to 1 when $n$ tends to infinity. So 
by the definition of an exposing functional, $f_n$ tends to $x_0^*$. 
By the same lemma $f_n(h_n)$ tends to~$1$, so $x_0^*(h_n)$ also tends
to~$1$. Hence every weak limit point of the sequence $(h_n)$ belongs 
to the  intersection of $A$ and  $S_{x_0^*}$, so this intersection 
is nonempty.
\end{proof}

\begin{prop} \label{expo2}
Let $E$ be a reflexive space.
\begsta
\item
If $E$ is USD-nonfriendly with USD-parameter ${<\alpha}$,
then there exists a functional  $x^* \in S(E^*)$ such
that for every strongly exposed  point $x_0^*$ of $B(E^*)$ 
the numerical
set $x^*(S_{x_0^*})$ contains the interval $[-1+\alpha, \allowbreak
1-\alpha]$.
\item
If $E$ is not USD-nonfriendly, then for every strongly
exposed  point $x_0^*$ of $B(E^*)$ the set $S_{x_0^*}$ has
diameter~$2$. Moreover,
for every $\delta > 0$ there exists a functional  $x^* \in S(E^*)$ such
that for every strongly exposed  point $x_0^*$ of $B(E^*)$ the numerical
set $x^*(S_{x_0^*})$ contains the interval $[-1+\delta, \allowbreak
1-\delta]$.
\endsta
\end{prop}

\begin{proof}
(a)
Let $A$ be a closed absolutely convex set 
which intersects all the sets $D\in \calD(E) $, but does not contain
$\alpha B(E)$.
By the Hahn-Banach theorem there exists a functional  $x^* \in S(E^*)$ such
that $|x^*(a)| < \alpha$ for every $a \in A $. 
We fix $y\in S(E)$ with $x^*(y)=-1$. 

Let $x_{0}^*\in S(E^*)$ be a strongly exposed point of $B(E^*)$.
As before, we denote an exposing evaluation functional by $x_{0}$.
Now $A\cap D(x_{0},y,\eps)\neq\emptyset$ for all $\eps>0$.
By Proposition~\ref{expo} and the evident equality         
$D(x_0 ,0, \eps) - y = D(x_{0} ,y, \eps)$
this implies that  the set  $A + y$ intersects 
$S_{x_0^*}$.
If $z_{1}$ is an element of this intersection, we see that 
$x^*(z_{1}) < \alpha -1$.

Likewise, since $D(-x_{0},0,\eps)= -D(x_{0},0,\eps)$, we find some
$z_{2}\in (-A-y)\cap S_{x_{0}^*}$; hence $x^*(z_{2})> -\alpha +1$.
Therefore, $[-1+\alpha, 1-\alpha] \subset x^*(S_{x_{0}^*})$.

(b) The argument is the same as in (a).
\end{proof}

This proposition allows us to estimate the USD-parameter of
finite-di\-men\-sional spaces.

\begin{prop} \label{endldim}
If $E$ is $n$-dimensional, then its USD-parameter is ${\ge 2/n}$.
\end{prop}

\begin{proof}
Assume that $\dim(E)=n$ and that its USD-parameter is ${<2/n}$; then
this parameter is strictly smaller than some $\alpha<2/n$. Choose
$x^*$ as in Proposition~\ref{expo2} so that 
\begin{equation} \label{eq1000}
[-1+\alpha, 1-\alpha] \subset x^*(S_{x_{0}^*})
\end{equation}
for every strongly exposed functional $x_{0}^*\in S(E^*)$.

We now claim that in any $\eps$-neighbourhood of $x^*$ there is some
$y^*\in B(E^*)$ which can be represented as a convex combination of
${\le n}$ strongly exposed functionals. First of all, the convex hull
of the set $\stexp B(E^*)$ of  strongly exposed functionals
is norm-dense in $B(E^*)$; in fact, this is true of any bounded
closed convex set in a separable dual space \cite[p.~110]{BL1}. Hence
for some $\|y_{1}^*-x^*\|<\eps$, $\lambda'_{1}, \dots,
\lambda'_{r}\ge0 $ with $\sum_{k=1}^r \lambda'_{k}=1 $ and $x_{1}^*,
\dots, x_{r}^*\in \stexp B(E^*)$
$$
y_{1}^* = \sum_{k=1}^r \lambda_{k}' x_{k}^*.
$$
Let $C=\co  \{ x_{1}^*, \dots, x_{r}^* \}$ and let $y^*$ be the
point of intersection of the segment $[y_{1}^*, x^*]$ with the relative
boundary of $C$, i.e., $y^*= \tau x^* + (1-\tau) y_{1}^*$
with $\tau = \sup\{ t\in [0,1]\dopu t x^* + (1-t) y_{1}^* \in C\}$.
Let $F$ be the face of $C$ generated by $y^*$; then $F$ is a convex
set of dimension~${< n}$. Therefore an appeal to Carath\'eodory's
theorem shows that $y^*$ can be represented as a convex combination
of no more than $n$ extreme points of $F$. But $\ex F \subset \ex C
\subset \{x_{1}^*, \dots, x_{r}^*\} \subset \stexp B(E^*)$, and our
claim is established.

We apply the claim with some $\eps< 2/n-\alpha$ to obtain some convex
combination $y^*= \sum_{k=1}^n \lambda_{k} x_{k}^*$ of $n$ strongly
exposed functionals  such that $\|y^*-x^*\|<\eps$. 
One of the coefficients must be ${\ge 1/n}$, say $\lambda_{n}\ge
1/n$. Now if $x\in S_{x_{n}^*}$,
\bea
x^*(x) &\ge& 
x^*(y)-\eps =
\sum_{k=1}^{n-1} \lambda_{k} x_{k}^*(x) + \lambda_{n} -\eps \\
&\ge&
-\sum_{k=1}^{n-1} \lambda_{k} + \lambda_{n} = -1+2\lambda_{n} -\eps 
\ge -1+2/n -\eps.
\eea
By \kref{eq1000} we have $-1+\alpha \ge -1+2/n -\eps$ which contradicts our
choice of $\eps$.
\end{proof}

For $\ell_{\infty}^n$ we can say more, namely, its USD-parameter is
worst possible.

\begin{prop} \label{ell_infty_n}
The USD-parameter of $\ell_{\infty}^{n}$ is $2/n$.
\end{prop}

\begin{proof}
The argument of
Proposition~\ref{propc0}(a) implies in the setting of
$\ell_{\infty}^{n}$ rather than $c_{0}$       that the USD-parameter of
$\ell_{\infty}^{n}$ is ${\le 2/n}$,
and the converse estimate follows from Proposition~\ref{endldim}.
\end{proof}

\section{Strong Daugavet and narrow operators in spaces of 
vector-valued functions}
\label{sec3}

Let $E$ be a Banach space and $X$ be a subspace of the space of all
bounded
$E$-valued 
functions defined on a set $K$, equipped with the sup-norm.
It will be convenient to use the following notation:
A disjoint pair $(U,V)$ of subsets of  $K$ is said to be
{\em interpolating}\/ for $X$ if for every $f,g \in X$ with $\|f\|<1$
and $\|g \chi_V\|<1$ there exists $h \in B(X)$ such that
$h=f$ on $U$ and $h=g$ on $V$. 

For arbitrary $V \subset K$ denote by $X_V$ the subspace of 
all functions from $X$ vanishing on $V$.

\begin{prop}  \label{ernst}
Let $X$ be as above and let
$(U,V)$ be an interpolating pair for~$X$. Then for every $f \in X$
$$
\dist(f,X_V) \le \sup_{t\in V} \|f(t)\| .
$$
\end{prop}

\begin{proof}
By the definition
of an interpolating pair, for an arbitrary $\eps > 0$ there exists
an element $h \in X$, $\|h\|< \sup_{t\in V} \|f(t)\|  + \eps$, such that 
$h=0$ on $U$ and $h=f$ on $V$. Then the element  
$f-h$ belongs to $X_V$, so 
$$
\dist(f, X_V) \le \|f - (f-h)\| = \|h\| < \sup_{t\in V} \|f(t)\|
+\eps,
$$
which completes the proof.
\end{proof}

\begin{lemma} \label{prox}
Let $X\subset \ell_{\infty}(K,E)$, $U,V\subset K$, $f\in S(X_{V})$ and
$\eps>0$. Assume that $U \supset \{ t\in K\dopu \|f(t)\|>1-\eps\}$ and
that $(U,V)$ is an interpolating pair for $X$. If 
$T$ is a strong Daugavet operator on $X$
and $g\in B(X)$, there is a function $h\in X_{V}$, $\|h\|\le 2+\eps$, 
satisfying
$$
\|Th\|<\eps, \ \|(g+h)\chi_{U}\|<1+\eps \ \text{and} \ 
\|(f+g+h)\chi_{U}\|> 2-\eps.
$$
\end{lemma}

\begin{proof}
Before we enter the proof proper, we formulate a number of technical
assertions that are easy to verify and will be needed later.

\begin{sublemma}\label{sublemma1}
If $T$ is a strong Daugavet operator on a Banach space $X$, 
if $1-\eta <\|x\| < 1+\eta$  and $\|y\| < 1+\eta$, then there is some
$z\in X$  such that
\[
\|x+y+z\| > 2-3\eta,\ \|y+z\| < 1+2\eta, \ \|Tz\|<\eta.
\]
\end{sublemma}

\begin{proof}
Choose $x_{0}\in S(X)$ and $y_{0}\in B(X)$ such that
$\|x_{0}-x\|<\eta$, $\|y_{0}-y\|<\eta$ and pick by Lemma~\ref{lem1.2}
$z\in D(x_{0},y_{0},\eta)$ such that $\|Tz\|<\eta$; this $z$ clearly
works.
\end{proof}

\begin{sublemma}\label{sublemma2}
If $\|x\|<1+\eta$, $\|y\|<1+\eta$ and
$\|(x+y)/2\| > 1-\eta$ in a normed space, 
then $\|\lambda x + (1-\lambda)y\| > 1-3\eta$
whenever $0\le\lambda \le1$.
\end{sublemma}

\begin{proof}
Should $\|\lambda x + (1-\lambda) y\| \le 1-3\eta$ for some $0\le
\lambda \le 1/2$, then, since $\lambda_{1}x + (1-\lambda_{1})(
\lambda x + (1-\lambda)y ) = (x+y)/2$ for $\lambda_{1}= (\frac12 -
\lambda)/(1-\lambda) \in [0,1/2]$, 
$$
\Bigl\| \frac{x+y}2 \Bigr\| \le \lambda_{1} (1+\eta) +
(1-\lambda_{1}) (1-3\eta) 
= 1 - (3-4\lambda_{1})\eta \le 1-\eta.
$$
(The case $\lambda > 1/2$ is analogous.)
\end{proof}

\begin{sublemma}\label{sublemma3}
If $\|y\|<1+\eta$ and $\|x+Ny\|/(N+1) >
1-3\eta$ in a normed space, then $\|(x+y)/2\| > 1- (2N+1)\eta$.
\end{sublemma}

\begin{proof}
Should $\|(x+y)/2\| \le 1-(2N+1)\eta$, then
\bea
\Bigl\| \frac{x+Ny}{1+N} \Bigr\| &\le&
\frac2{1+N} \Bigl\| \frac{x+y}2 \Bigr\| + \Bigl( 1 - \frac2{1+N}
\Bigr) \|y\| \\
&\le&
\frac2{1+N} \bigl( 1- (2N+1)\eta \bigr) + \Bigl( 1 - \frac2{1+N}
\Bigr) (1+\eta) = 1-3\eta.
\eea
\end{proof}

To start the actual proof we may assume that $\|T\|=1$. Fix $N>6/\eps$
and $\delta > 0$ such that $2(2N+1)9^N\delta < \eps$; let
$\delta_{n}= 9^n \delta$ so that $(2N+1) \delta_{N} < \eps/2$.
Put $f_{1}=f$, $g_{1}=g$ and pick $h_{1}\in X$ such that
$$
\|f_{1}+g_{1}+h_{1}\| > 2-\delta_{1}, \ \|g_1+h_{1}\| <
1+2\delta_{0}, \ \|Th_{1}\|<\delta_{0}.
$$

We are going to 
construct functions $f_{n},g_{n},h_{n}\in X$ by induction so as to
satisfy
\begsta
\item
\label{prox1}
$f_{n+1} = \frac{1}{n+1}(f_1 + \sum_{k=1}^n(g_k + h_k)) =
\frac{n}{n+1} f_{n} + \frac1{n+1}(g_{n}+h_{n})$,
$1-3\delta_{n} < \|f_{n+1}\| <1+\delta_{n}$,
\item 
\label{prox2} 
 $g_{n+1}= g_1$ 
on $U$   and $g_{n+1}=g_n +  h_n $ $( = g_{1} +h_1+ \dots  + h_n)$ on
$V$, $\|g_{n+1}\| < 1+\delta_{n}$,
\item 
\label{prox3} 
$\| f_{n+1} + g_{n+1} +h_{n+1} \| > 2-\delta_{n+1}$,
$1-2\delta_{n} < \| g_{n+1} + h_{n+1} \| < 1+6\delta_{n} <
1+\delta_{n+1}$, $\|Th_{n+1}\| < 3 \delta_{n}$.
\endsta

Suppose that these functions have already
been constructed for the indices
$1,\dots,n$. We then define $f_{n+1}$ as in (a). Since by induction
hypothesis $\|f_{n}\|<1+\delta_{n-1}$ and $\|g_{n}+h_{n} \| <
1+\delta_{n}$, we clearly have $\|f_{n+1}\| < 1+\delta_{n}$. From
$\|f_{n}+g_{n}+h_{n}\| > 2-\delta_{n}$ we conclude using
Sublemma~\ref{sublemma2} (with $\eta=\delta_{n}$) that $\|f_{n+1}\|>
1-3\delta_{n}$. Thus (a) is achieved. To achieve (b) it is enough to
use that $(U,V)$ is interpolating along with the induction hypothesis
that $\|g_{n}+h_{n}\| < 1+\delta_{n}$. Finally (c) follows from
Sublemma~\ref{sublemma1} with $\eta= 3\delta_{n}$.

Next we argue that
$$
\biggl\|f_1 + \frac{1}{N}\sum_{k=1}^N (g_k + h_k) \biggr\| > 2-\eps/2.
$$
This follows from Sublemma~\ref{sublemma3}, (c) and (a) and our
choice of $\delta$.
But for $t\notin U$ we can estimate
 $$
\biggl\|f_1(t) + \frac{1}{N}\sum_{k=1}^N (g_k (t) + h_k (t)) \biggr\| 
\le 1-\eps + 1-\delta_{N} \le 2-2\eps,
$$
therefore, letting $w= \frac1N \sum_{k=1}^N h_{k}$,
$$
\| (f+g+w)\chi_{U} \| = 
\biggl\| \biggl( f_1 + \frac{1}{N}\sum_{k=1}^N (g_k + h_k) \chi_{U}
\biggr) \biggr\| > 2-\eps/2.
$$
Furthermore we have the estimates
\bea
 \| (g+w)\chi_{U} \| &=&
\biggl\| \frac1N \sum_{k=1}^N (g_{k}+h_{k}) \chi_{U} \biggr\| \le
1+\delta_{N} < 1+\eps/2 ,\\
\|Tw\| &\le& 
\frac1N \sum_{k=1}^N \|Th_{k}\| <
3\delta_{N-1} =  \frac13 \delta_{N} < \eps/2, \\
 \|h_{k}\| &\le&
 \|g_{k}+h_{k}\| + \|g_{k}\| \le 2+2\delta_{k} \le 2+2\delta_{N}        \le
 2+\eps/2 ,\\
\|w\| &\le& 
\frac1N \sum_{k=1}^N \|h_{k}\| \le 
2+\eps/2 , 
\eea
and for $t\in V$
$$
\|w(t)\| 
= \frac1N \|g_{N+1}(t) - g_{1}(t)\| \le \frac{2+\delta_{N}}N 
< \frac 3N < \eps/2.
$$
By Proposition~\ref{ernst} and the above we see that 
$\dist(w,X_{V}) < \eps/2$. Hence it is left 
to replace $w$ by an element $h\in X_{V}$,
$\|h-w\|\le \eps/2$, to finish the proof. 
\end{proof}

Let us remark that the conditions of Lemma~\ref{prox} are
fulfilled for an arbitrary compact Hausdorff space $K$,
for a closed subset $V \subset K$ and for $X=C(K,E)$ as well as for
$X=C_w(K,E)$. Here is another example.

\begin{cor}\label{cor3.2a}
If $X= X_{1} \oplus_{\infty} X_{2}$ and $T\in \SD(X)$, then
$\rest{T}{X_{1}} \in \SD(X_{1})$.
\end{cor}

To see this let $K= \ex B(X^*)$, $K_{1}= \ex B(X_{1}^*)$, $K_{2}= \ex
B(X_{2}^*)$ so that $K=K_{1} \cup K_{2}$ and $X\subset
\ell_{\infty}(K)$ canonically. It is left to apply Lemma~\ref{prox}
with the interpolating pair $(K_{1},K_{2})$.
A direct proof of
Corollary~\ref{cor3.2a} is given in \cite{BKSW}.

\begin{thm} \label{cnarvec}
Let $K$ be a compact Hausdorff space,
$E$ a Banach space and  $T$ 
an operator on $X=C(K,E)$. Then the following conditions are equivalent:
\begin{enumerate}
\item \label{cnar2}
 $T \in \SD(X)$.
\item \label{cnar3}
For every closed subset $V \subset K$, every $x \in S(E)$, every 
$y \in B(E)$  and  every $\eps >0$ there exists an open subset 
$W \subset K \setminus V$, an element $e \in E$ with $\|e + y\| < 1 + \eps$,  
$\|e + y + x\| > 2- \eps$, and a  function $h \in X_V$, $\|h\|\le
2+\eps$, such that $\|Th\| < \eps$ and $\|e-h(t)\|< \eps$ for $t \in W$. 
\item \label{cnar4}
For every closed subset $V \subset K$, every $x \in S(E)$, every 
$y \in B(E)$  and  every $\eps >0$ there exists a function $f \in X_V$
such that $\|Tf\| < \eps$, $\|f + y\| < 1 + \eps$,  
$\|f + y + x\| > 2- \eps$.
\end{enumerate}
If $K$ has no isolated points, then these conditions are equivalent
to
\begin{enumerate}  \setcounter{enumi}{3}  
\item \label{cnar1}
 $T \in \Nar(X)$.
\end{enumerate}
\end{thm}

\begin{proof}
The implication~\kref{cnar2} $\Rightarrow $ \kref{cnar3} follows 
from Lemma~\ref{prox}, as follows.
 Let us apply 
Lemma~\ref{prox} to  ${\eps}/{4} >0$, 
$g=\chi_{K} \otimes y$, $f=f_1 \otimes x \in S(X)$, 
where $f_1$ is a positive scalar function vanishing on $V$, and 
$U=\{t\in K\dopu \|f(t)\| > 1-\eps/4 \}$. Then for 
$h \in X_V$ which we get from Lemma~\ref{prox} let us find a point
$t_{0} \in U$ such that 
$\|(f+g+h)(t_{0})\|= \|(f+h)(t_{0}) + y \|> 2 -{\eps}/{4}$.
Because $\|h(t_{0}) + y \| < 1 + {\eps}/{4}$ we have 
$\|f(t_{0}) \|> 1 -{\eps}/{2}$, i.e., $\|f(t_{0})- x \| < \eps/{2}$.
Now select an open neighbourhood $W \subset U$
of $t_{0}$ such that $\|f(\tau)- x \| < {\eps}/{2}$ for all $\tau \in W$
and put $e = h(t_{0})$.

To prove the implication~\kref{cnar3} $\Rightarrow $ \kref{cnar4}
let us fix a positive $\eps < {1}/{10}$, $\delta < {\eps}/{4}$
and $N > 6 + 2/\eps$. 
Now apply inductively condition~\kref{cnar3} to obtain
elements $x_k, y_k, e_k$, $x_1=x$, $ y_k=y$, $k =1, \dots  , N$, open subsets 
$W_1 \supset W_2 \supset \dots $, closed subsets 
$V_{k+1} = K \setminus W_k$, $V_1 =V$ and functions $h_k \in X_{V_k}$ with
the following properties:

\begsta
\item \label{cnar1.1}
$\displaystyle
x_{n+1} = \frac{x + \sum_{k=1}^n(y_k + e_k)}{\|x + \sum_{k=1}^n(y_k + e_k)\|}
\in S(E)$, 
\item \label{cnar1.2}
$\|e_k + y_k\| < 1 + \delta$,  
$\|e_k + y_k + x_k\| > 2- \delta$,
\item \label{cnar1.3}
$h_k \in X_{V_k}$, $\|h_k(t) - e_k\| <  {\eps}/{4}$ for all 
$t \in W_k$, $\|h_k\| \le 2 + \eps$, and $\|Th_k\| <  \eps$.
\endsta

By an argument similar to the one in Lemma~\ref{prox}, we have for a
proper choice of $\delta$
$$
\biggl\|x + y + \frac{1}{N} \sum_{k=1}^N e_k \biggr\| = \biggl\|x +
\frac{1}{N}\sum_{k=1}^N(y_k + e_k) \biggr\| > 2 - \frac{\eps}{2}.
$$
Let us put $f=\frac{1}{N}\sum_{k=1}^N h_k$. 
Then the last inequality  and (c) of 
our construction yield that $f \in X_V$, 
$\|f + y + x\| > 2- \eps$ and $\|Tf\| < \eps$. The only thing left
to do now is to estimate $\|f + y\|$ from above.
If $t \in V$, then $\|f(t) + y\|=\|y\| \le 1$. If 
$t \in W_n \setminus W_{n+1}$ for some $n$ then 
$$
\|f(t) + y\|= \biggl\| \frac{1}{N}\sum_{k=1}^n h_k(t) + y \biggr\|= 
\biggl\|\frac{1}{N}\sum_{k=1}^n (h_k(t) + y) \biggr\|
$$
In this sum all the summands except for the last one satisfy the 
inequality $\|h_k(t) + y\| \le 1 + {\eps}/{2}$ and the last summand
$h_n(t) + y$ is bounded by $3 + \eps$. So
$$
\|f(t) + y\| \le \frac{1}{N}\sum_{k=1}^{n-1} 
\Bigl(1 + \frac{\eps}{2}\Bigr) + \frac{1}{N}(3 +
\eps) \le  1 + \frac{\eps}{2} + \frac{1}{N}(3 + \eps) \le  1 + \eps.
$$
The same estimate holds for $t \in W_N$. 

To prove the implication \kref{cnar4} $\Rightarrow$ \kref{cnar2} fix
$f,g\in S(X)$ and $0<\eps <1/10$. Pick a point $t\in K$ with
$\|f(t)\| > 1-\eps/4$ and a neighbourhood $U$ of $t$  such that
$$   
\|f(t)- f(\tau)\| + \|g(t) - g(\tau)\| < \frac\eps4 \qquad\forall
\tau\in U.
$$   
Denote $x= f(t)/\|f(t)\|$ and $y=g(t)$ and apply
condition~\kref{cnar4} to obtain a function $h\in X_{V}$ such that
$\|Th\|<\eps$, $\|h+y\| <1+\eps/4$ and $\|h+y+x\| >2-\eps/4$. For
this $h$ we have $\|h+g\|<1+\eps$ and $\|h+g+f\|>2-\eps$, so $T\in
\SD(X)$.

Let us now pass to the case of a perfect compact $K$. 
The implication~\kref{cnar1} $\Rightarrow $ \kref{cnar2} is evident.

The proof of the remaining implication~\kref{cnar4} $\Rightarrow 
$ \kref{cnar1} is similar to that of \kref{cnar4} $\Rightarrow 
$ \kref{cnar2}. 
Namely, let  $f, g \in S(X)$, $x^* \in X^*$ and let $\eps > 0$ be small. 
We have to show that there is an element $h \in X$ 
such that
\begin{equation}
\label{eq006} \|f + g + h\|> 2 - \eps,  \quad \| g + h\|< 1+ \eps 
\end{equation}
and 
\begin{equation}
\label{eq009} \|Th\|+|x^*h| < \eps. 
\end{equation}
To this end let us pick a closed subset 
$V \subset K$ (whose complement $K \setminus V$  we denote by $U$)  
and a point $t \in U$ in such a way that
$\|f(t)\| > 1 - \eps/{4}$,  
\begin{equation}
\label{eq007} |x^*|_{X_V} < \frac{\eps}{4}, 
\end{equation}
and  for every $\tau \in U$ 
\begin{equation}
\label{eq008} \|f(t) - f(\tau)\| + \|g(t) - g(\tau)\| < \frac{\eps}{4}.
\end{equation}
Denote $x={f(t)}/{\|f(t)\|}$, $y=g(t)$ and apply condition~\kref{cnar4} 
to obtain a function $h \in X_V$ such that 
$\|Th\| <{\eps}/{4} $, $\|h + y\| < 1 + {\eps}/{4}$  and
$\|h + y + x\| > 2- {\eps}/{4}$. For this $h$ 
(\ref{eq006}) follows from (\ref{eq008}) and (\ref{eq009}) 
follows from (\ref{eq007}).
\end{proof}

In \cite{KadSW2} we have introduced the tilde-sum of two operators
$T_{1}\dopu X\to Y_{1}$, $T_{2}\dopu X\to Y_{2}$ by
$$
T_{1} \mytilde+ T_{2} \dopu X \to Y_{1} \oplus_{1} Y_{2},
\ x\mapsto (T_{1}x, T_{2}x).
$$
There we proved that the $\mytilde+$-sum and therefore also the
ordinary sum of two narrow operators on $C(K)$ is narrow (another
proof will be given in the next section), and we inquired whether
this is so on any space with the \DP. We are now in a position to
provide a counterexample.

Let $T\dopu E\to F$ be an operator on a Banach space. 
By  $T^K$ let us denote
the corresponding ``multiplication'' or ``diagonal''
operator $T^K\dopu C(K,E) \to C(K,F)$ defined by
$$
(T^K f)(t) = T(f(t)).
$$

\begin{prop} \label{mult1}
$T^K \in \SD(C(K,E))$ if and only if  $T\in \SD(E)$.
\end{prop}

\begin{proof}
Criterion~\kref{cnar4} of Theorem~\ref{cnarvec} immediately provides
the proof.
\end{proof}

Here is the announced counterexample.

\begin{thm}\label{theo6.5}
There exists a Banach space $X$ for which $\Nar(X)$ does not form a
semigroup under the operation $\mytilde{+}$; in fact,
$C([0,1], \ell_{1})$ is such a space.
\end{thm}

\begin{proof}
The key feature of $\ell_{1}$   is that 
$\SD(\ell_{1})$ is not a $\mytilde +$-semigroup; for we have shown in
Proposition~\ref{propc0}(b)
that $x_{1}^*(x) = \sum_{n= 1}^\infty x(n)$ and 
$x_{2}^*(x) = x(1) - \sum_{n= 2}^\infty x(n)$ define strong Daugavet
functionals on $\ell_{1}$, but $x_{1}^* + x_{2}^*\dopu x \mapsto
2x(1)$ is not in $\SD(\ell_{1})$ and hence $x_{1}^* \mytilde+ x_{2}^*$
is not, either. 

Now if $\SD(E)$ fails to be a $\mytilde{+}$-semigroup, pick
$T_1,T_2 \in \SD(E)$ with
$T_1 \mytilde{+} T_2 \notin \SD(E)$.
Put $X=C(K, E)$ for a perfect compact Hausdorff space $K$; 
then by Proposition~\ref{mult1} and
Theorem~\ref{cnarvec} 
$T_1^{K},T_2^{K} \in \Nar(X))$, but 
$T_1^{K} \mytilde{+} T_2^{K} \notin \Nar(X)$.
\end{proof}

Another example of a space for which $\SD(E)$ is no
$\mytilde+$-semigroup is  $E=L_{1}[0,1]$. This  is much subtler than
for $\ell_{1}$; the proof 
is presented in \cite[Th.~6.3]{KadSW2}. This example has the
additional benefit of involving a space with the \DP; by
Theorem~\ref{theo6.5}, however, $E=C([0,1],\ell_{1})$ is another
example of this kind.

\section{Narrow and $C$-narrow operators on  $C(K,E)$}
\label{sec4}

The following definition extends the notion of a $C$-narrow operator 
studied in \cite{KadPop} and \cite{KadSW2} to the vector-valued
setting.

\begin{definition}
\label{narrow}
An operator $T\in L(C(K,E),W)$ is called 
\textit{$C$-narrow} if there is a constant $\lambda$ such that
given any $\varepsilon >0$, $x\in S(E)$ and open set $U\subset K$ there is a
function $f\in C(K,E)$, $\Vert f\Vert \le \lambda$, satisfying the following
conditions:
\begsta
\item
$\supp(f)\subset U$,
\item
$f^{-1}(B(x,\varepsilon ))\neq \emptyset $, where $
B(x,\varepsilon )= \{z\in E\dopu \|z-x\|<\eps\}$,
\item
$\Vert Tf\Vert <\varepsilon $.
\endsta
\end{definition}

As the following proposition shows, condition~(b) of the previous definition
can be substantially strengthened. In particular, the size of the
constant $\lambda$ is immaterial; but introducing this constant in
the definition allows for more flexibility in applications. Also,
Proposition~\ref{reduction} shows that for $E=\R$ the new notion of
$C$-narrowness coincides with the one from \cite{KadSW2}. 

\begin{prop}
\label{reduction}
If $T$ is a $C$-narrow operator, then 
for every $\eps > 0$, $x\in S(E)$ and open set $U\subset K$ 
there is a function $f$ 
of the form $g\otimes x$, where $g\in C(K)$, $\supp(g)\subset U$,
$\Vert g\Vert =1$ and $g$ is
nonnegative, such that $\Vert Tf\Vert <\varepsilon $.
\end{prop}

\begin{proof}
Let us fix $\varepsilon >0$, an open set $U$ in $K$ and 
$x\in S(E)$. By Definition~\ref{narrow} we find a function $f_{1}\in C(K,E)$
as described there
corresponding to $\varepsilon $, $U$ and $x$. Put $U_1=U$ and $%
U_2=f_1^{-1}(B(x,\frac 12))$. As above, there is a function $f_2$
corresponding to $\eps$, $U_2$ and~$x$. 
We denote $U_3=f_2^{-1}(B(x,\frac 14))$ and continue
the process.
In the $r^{\rm th}$ step we get the set $U_r=f_{r-1}^{-1}(B(x,\frac
1{2^{r-1}}))$
and apply Definition~\ref{narrow} to obtain a function $f_r$ corresponding
to $U_r$.

Choose $n\in \mathbb{N}$ so that $(\lambda+ 2)/n<\varepsilon $
and put $f=\frac
1n(f_1+f_2+\dots +f_n)$. Now using the Urysohn Lemma we find a
continuous function $g$
satisfying $\frac{k-1}n\leq g(t)\leq \frac kn$ for all $t\in U_k$, $%
k=1,\dots ,n$, $\|g\|=1$
and vanishing outside $U_1$. We claim that $\Vert
f-g\otimes x\Vert <\varepsilon $. Indeed, by our construction, 
if $t\in K\setminus U_{1}$, then $\Vert (f-g\otimes x)(t)\Vert
=0$, and if $t\in U_{k} \setminus U_{k+1}$ (with the understanding
that $U_{n+1}$ stands for $\emptyset$), then
\begin{eqnarray*}
\Vert (f-g\otimes x)(t)\Vert 
 &=&
\Bigl\| \frac 1n(f_1+\dots +f_k)(t)-g(t)\cdot x\Bigr\| \\
&\leq &
\Bigl\| \frac1n \bigl( (f_{1}(t)-x) + \dots + (f_{k-1}(t)-x) +
f_{k}(t) \bigr) \Bigr\| + \frac1n \\
&\le&
\frac1n\Bigl(  \frac 12 +\dots +\frac 1{2^{k-1}}+ \lambda \Bigr) +
\frac1n  < \frac {\lambda+2}n < \varepsilon .
\end{eqnarray*}
Moreover, 
$$
\Vert Tf\Vert \leq \frac 1n\left( \Vert Tf_1\Vert +\Vert
Tf_2\Vert +\dots +\Vert Tf_n\Vert \right) <\varepsilon .
$$ 
Thus $\Vert
T(g\otimes x)\Vert <\varepsilon +\varepsilon \Vert T\Vert $, and since $%
\varepsilon $ was chosen arbitrarily, we are done. 
\end{proof}

Another way to express this proposition is to say that $T\dopu
C(K,E)\to W$ is $C$-narrow if and only if, for each $x\in E$, the
restriction $T_{x}\dopu C(K)\to W$, $T_{x}(g)= T(g\otimes x)$, is
$C$-narrow.

\begin{prop} \label{Cvs.nar} \mbox{}
\begsta
\item
Every $C$-narrow operator on  $C(K,E)$ is a strong Daugavet
operator. Hence, in the case of a perfect compact $K$ every
$C$-narrow operator on  $C(K,E)$ is narrow.
\item
If $E$ is a separable USD-nonfriendly space, 
then every strong Daugavet operator on $C(K,E)$ is $C$-narrow.
\item
If  every strong Daugavet operator on
$C(K,E)$ is $C$-narrow, then 
$E$ is SD-nonfriendly.
\endsta
\end{prop}

\begin{proof}
(a) Let $T$ be $C$-narrow.
We will use criterion (\ref{cnar4})  of Theorem~\ref{cnarvec}. 
Let $F \subset K$ be a closed subset, 
 $x \in S(E)$,  
$y \in B(E)$  and  $\eps >0$. According to Proposition~\ref{reduction}  
there exists a function $f $ vanishing on $F$ of the
form $g \otimes (x-y)$, where $g \in C(K)$, $\Vert g\Vert =1$ and $g$ is
nonnegative, such that $\|Tf\| < \eps$. Evidently this $f$ 
satisfies all the demands of criterion (\ref{cnar4}) in
Theorem~\ref{cnarvec}.

(b) Let $T$ be a strong Daugavet operator, and suppose $E$ is
separable. Let $U\subset K$ be a nonvoid open subset.
Given $x,y\in S(E)$ and $\eps'>0$ we define
\bea
O(x,y,\eps') &=&
\{t\in U \dopu \ \exists f\in C(K,E)\dopu \supp f \subset U , \
\|f+y\|<1+\eps', \\
&& \qquad \qquad 
\|f(t)+y+x\|> 2-\eps', \ \|Tf\|<\eps' \}.
\eea 
This is an open subset of $K$, and by
Theorem~\ref{cnarvec}\kref{cnar4} it is dense in $U$. Now pick a countable
dense subset $\{(x_{n},y_{n})\dopu n\in\N\}$ of $S(E)\times S(E)$ and
a null sequence $(\eps_{n})$. Then by Baire's theorem, $G:=
\bigcap_{n} O(x_{n},y_{n},\eps_{n})$ is nonempty.

Let   $\eps >0$, and fix $t_{0}\in
G$. We  denote by $A(U,\eps)$ the
closure of 
$$
\{f(t_{0})\dopu f\in C(K,E),\ \|f\|<2+\eps,\ \|Tf\|<\eps, \ \supp
f\subset U \};
$$
this is an absolutely convex set. We claim that $A(U,\eps)$
intersects each set  $D(x,y,\eps')\in {\cal D}(E)$. 
Indeed, if $\|x_{n}-x\|<\eps'/4$, $\|y_{n}-y\|<\eps'/4$,
$\eps_{n}<\eps'/2$ and $\eps_n<\eps$, then
for a function $f_{n}$ as appearing in the definition of
$O(x_{n},y_{n},\eps_{n})$ we have 
$f_{n}(t_{0}) \in A(U,\eps) \cap D(x_{n},y_{n},\eps_{n}) \subset 
A(U,\eps) \cap D(x,y,\eps') $. 

Since $E$ is USD-nonfriendly, say with parameter $\alpha$, the set
$A(U,\eps)$ contains $\alpha B(E)$.
This implies that 
$T$ satisfies the definition of a $C$-narrow operator with constant
$\lambda = 3/\alpha$.

(c) Let $T\in \SD(E)$; then by Proposition~\ref{mult1} $T^K$ is a
strong Daugavet operator on $C(K,E)$. But
$$
\bigl( T^K (g\otimes e) \bigr) (t) = T \bigl( (g\otimes e)(t) \bigr)
= g(t) Te,
$$
hence $T^K$ is not $C$-narrow unless $T=0$.
\end{proof} 

The example $E=c_{0}$ shows that the converse of (b) is false. We
have already pointed out in Proposition~\ref{propc0}(a) that $c_{0}$
fails to be USD-nonfriendly; yet every \sDo\ on $C(K,c_{0})$ is
$C$-narrow. To see this we first remark that it is enough to check
the condition spelt out in Proposition~\ref{reduction} for $x$ in a
dense subset of $S(E)$. In our context we may therefore assume that the
sequence $x$ vanishes eventually, say $x(n)=0$ for $n>N$.  If we
write $c_{0}= \ell_{\infty}^N \oplus_{\infty} Z$, with $Z$ the space
of null sequences supported on $\{N+1, N+2, \dots\}$, we also have
$C(K, c_{0}) = C(K, \ell_{\infty}^N) \oplus_{\infty} C(K,Z)$. By
Corollary~\ref{cor3.2a}
the restriction of any \sDo\ $T$ on $C(K, c_{0}) $ to
$C(K, \ell_{\infty}^N) $ is again a \sDo, and hence it is $C$-narrow,
for $\ell_{\infty}^N$ is USD-nonfriendly (Proposition~\ref{prop2.7}).
This implies that $T$ is $C$-narrow.

We do not know whether (c) is actually an equivalence.

One of the fundamental 
properties of $C$-narrow operators is stated in our next
theorem. 

\begin{thm}
\label{sum}
Suppose that operators $T$, $T_n\in L(C(K,E),W)$ 
are such that the series $\sum_{n=1}^\infty w^{*}(T_nf)$
converges absolutely  to $w^{*}(Tf)$, for every 
$w^{*}\in W^{*}$ and $f\in C(K,E)$. If all the $T_n$ are $C$-narrow,
then so is $T$. In particular, the sum of two
$C$-narrow operators is a $C$-narrow operator.
\end{thm}

\begin{cor}\label{cor7.5}
A pointwise unconditionally convergent sum of narrow operators on
$C(K,E)$ is a narrow operator itself if $E$ is separable and USD-nonfriendly.
\end{cor}

Indeed, this follows from Theorem~\ref{sum} and
Proposition~\ref{Cvs.nar}. We remark that the case of a sum of two
narrow operators on $C(K)$ was treated earlier in \cite{KadPop} and
\cite{KadSW2}, but the assertion about infinite sums is new even there.
It was proved in \cite{KadSSW} for a pointwise unconditionally
convergent sum $T=\sum_{n=1}^\infty T_{n}$  on a space with the \DP\ that
$$
\|\Id+T\|\ge1
$$      
whenever $\|\Id+S\|=1+\|S\|$ for every $S$ in the
linear span of the $T_{n}$. In the context of Theorem~\ref{sum}
we even obtain 
\begin{equation}   \label{eq6.3}
\|\Id+T\| = 1+\|T\|
\end{equation}
when all the $T_{n}$ are narrow on $C(K)$. In particular,
the identity on $C(K)$ cannot be represented as
an unconditional sum of narrow operators, since obviously \kref{eq6.3}
fails for $T=-\Id$. This last consequence shows for an unconditional
Schauder decomposition $C(K)= X_{1} \oplus X_{2} \oplus \dots$ with
corresponding projections $P_{1}, P_{2}, \dots$ that one of the $P_{n}$
must be non-narrow, since $\Id = \sum_{n=1}^\infty P_{n}$ pointwise
unconditionally. Hence one of the $X_{n}$ must be
infinite-dimensional if $K$ is a perfect compact Hausdorff space. In
fact, one of the $X_{n}$ must contain a copy of $C[0,1]$ and
therefore be isomorphic to $C[0,1]$ by a theorem due to
Pe{\l}czy\'nski \cite{Pel-CSsub} if $K$ is in addition metrisable; 
 see \cite{KadPop} and \cite{KadSSW} for more results along these lines.

We now turn to the proof of Theorem~\ref{sum} for which we need
an auxiliary concept. A similar idea has appeared in \cite{KadPop}.

\begin{definition}
\label{van}
Let $G$ be a closed $G_\delta $-set in $K$ and $T\in L
(C(K),W)$. We say that $G$ is a \textit{vanishing set} 
of $T$ if there is a sequence
of open sets $(U_i)_{i\in \mathbb{N}}$ in $K$ and a sequence of functions $%
(f_i)_{i\in \mathbb{N}}$ in $S(C(K))$ such that
\begsta
\item
$G={\bigcap_{i=1}^\infty}U_i$;
\item
$\supp(f_i)\subset U_i$;
\item
${\lim_{i\rightarrow \infty } }f_i=\chi _G$ pointwise;
\item
${\lim_{i\rightarrow \infty } }\Vert Tf_i\Vert =0$.
\endsta
The collection of all vanishing sets of $T$ is denoted by $\van T$.
\end{definition}

Let $T\in L(C(K),W)$. By the Riesz Representation Theorem, $%
T^{*}w^{*}$ can be viewed as a regular measure on the Borel subsets of $K$
whenever $w^{*}\in W^{*}$. For convenience, we denote it by $T^{*}w^{*}$
as well.

\begin{lemma}
\label{aux1} 
Suppose $G$ is a closed $G_\delta $-set in $K$ and $T\in 
L(C(K),W)$. Then $G\in \van T$ if and only if $%
T^{*}w^{*}(G)=0$ for all $w^{*}\in W^{*}$.
\end{lemma}

\begin{proof} 
Let $G\in \van T$, and pick functions $(f_i)_{i\in \mathbb{N}}$ as in
Definition~\ref{van}. Then by the Lebesgue Dominated Convergence Theorem,
for any given $w^{*}\in W^{*}$ we have 
$$
T^{*}w^{*}(G) = \int_K  \! \chi _G \, dT^{*}w^{*}= \lim_{i\rightarrow
\infty } \int_K \! f_i \, dT^{*}w^{*} 
= \lim_{i\rightarrow \infty } w^{*}(Tf_i)=0.
$$

Conversely,
let $(U_i)_{i\in \mathbb{N}}$ be a sequence of open sets in $K$
such that $\overline{U}_{i+1}\subset U_i$ and $G=
\bigcap_{i=1}^\infty U_i$. By the Urysohn Lemma there exist functions $%
(f_i)_{i\in \mathbb{N}}$ having the following properties: $0\leq f_i(t)\leq 1$
for all $t\in K$, $\supp (f_i)\subset U_i$, and $f_i(t)=1$ if $t\in 
\overline{U}_{i+1}$. Clearly, $\lim_{i\rightarrow \infty } 
f_i=\chi _G$ pointwise and 
\[
\lim_{i\rightarrow \infty } w^{*}(Tf_i) =
\lim_{i\rightarrow \infty } T^{*}w^{*}(f_i) =
T^{*}w^{*}(G)=0
\]
whenever $w^{*}\in W^{*}$. This means that the sequence $(Tf_i)_{i\in \N%
}$ is weakly null. Applying the Mazur Theorem we finally obtain a sequence
of convex combinations of the functions $(f_i)_{i\in \N}$ which satisfies
all the conditions of Definition~\ref{van}.

This completes the proof. 
\end{proof}

\begin{lemma}
\label{aux2}
An operator $T\in L(C(K),W)$ is $C$-narrow if and only if
every nonvoid open set $U\subset K$ contains a nonvoid vanishing set of $T$.
Moreover, if 
$(T_n)_{n\in \N}\subset L(C(K),W)$ is a sequence of $C$-narrow
operators, every open set $U\neq\emptyset $ contains a set
$G\neq\emptyset $ that is
simultaneously a vanishing set for all $T_n$. 
\end{lemma}

\begin{proof} 
We first prove the more general ``moreover'' part.
Put $U_{1,1}=U$. By the definition of a $C$-narrow operator and
Proposition~\ref{reduction} there is a function 
$f_{1,1}\subset S(C(K))$ with $\supp (f_{1,1})\subset U_{1,1}$, 
$U_{1,2}:=f_{1,1}^{-1}(\frac 12,1]\neq \emptyset$ and 
$\Vert T_1f_{1,1} \Vert <\frac 12$. Obviously,   
$\overline{U}_{1,2}\subset f_{1,1}^{-1}[\frac12,1]\subset U_{1,1}$.  
Again applying the definition we find $f_{1,2}\in S(C(K))$ with 
$\supp (f_{1,2})\subset U_{1,2}$, 
$ U_{2,1}=f_{1,2}^{-1}(\frac 23,1]\neq \emptyset $ and 
$\Vert T_1f_{1,2}\Vert <\frac 13$. As above 
$\overline{U}_{2,1}\subset U_{1,2}$.

In view of the $C$-narrowness of $T_2$ there exists a function  
$f_{2,1}\in S(C(K))$ with $\supp (f_{2,1})\subset U_{2,1}$, $%
U_{1,3}=f_{2,1}^{-1}(\frac 23,1]\neq \emptyset $ and $\Vert T_2f_{2,1}\Vert
<\frac 13$. In the next step we construct $f_{1,3}\in S(C(K))$ such that 
$U_{2,2}=f_{1,3}^{-1}(\frac 34,1]\neq \emptyset $ and $\Vert T_1f_{1,3}\Vert
<\frac 14$.

Proceeding in the same way, in the $n^{\rm th}$ 
step we find a set of functions 
$(f_{k,l})_{k+l=n}\subset S(C(K))$ and nonempty open sets $(U_{k,l})_{k+l=n}$
in $K$ such that $\supp (f_{k,l})\subset U_{k,l}$, $\Vert
T_kf_{k,n-k}\Vert <\frac 1n$ and $U_{k,l}=f_{k-1,l+1}^{-1}(\frac{n-1}n,1]$,
if $k\neq 1$. Then we put $U_{1,n}=f_{n-1,1}^{-1}(\frac{n-1}n,1]$ to start
the next step.

It remains to show that the set $G= \bigcap_{k,l\in \N} U_{k,l}=
\bigcap_{k,l\in \N} \overline{U}_{k,l}$
is as desired. Indeed, $G$ is clearly a nonempty
closed $G_\delta $-set and $G=%
\bigcap_{i=1}^\infty U_{n,i}$ for every $n\in \N$. 
It is easily seen that the sequences $(f_{n,i})_{i\in \N}$ and $%
(U_{n,i})_{i\in \N}$ meet the conditions of Definition~\ref{van} for
the operator $T_n$. So, $G\in \van T_n$ for every $n\in \N$.

To prove the converse,
let $U\neq\emptyset $ be any open set in $K$ and let $\varepsilon >0$. By 
assumption on $\van T$
we can find a closed $G_\delta $-set $\emptyset \neq
G\subset U$, $G\in \van T$. 
Consider the open sets $(U_i)_{i\in \N}$ and functions $%
(f_{i})_{i\in \N}$ provided by Definition~\ref{van}. For sufficiently
large $i\in \N$  we have $U_i\subset U$ and $\Vert Tf_i\Vert
<\varepsilon $ so that $f_i$ may serve as a function as 
required in Definition~\ref{narrow}.

This finishes the proof. 
\end{proof}

Now we are in a position to prove Theorem~\ref{sum}.

\smallskip\noindent
\emph{Proof of Theorem~\ref{sum}.} 
By virtue of Proposition~\ref{reduction}, 
we may assume that $E=\R$. By Lemma~\ref{aux2} it suffices
to show that $\bigcap_{n=1}^\infty \van T_n\subset \van T$.

Suppose $G\in \bigcap_{n=1}^\infty \van T_n$.
According to Lemma~\ref{aux1} we need to prove that $T^{*}w^{*}(G)=0$ for
all $w^{*}\in W^{*}$. By the condition of the theorem, 
the series ${\sum_{n=1}^\infty}T_n^{*}w^{*}$ is weak$^{*}$-unconditionally
Cauchy and hence weakly unconditionally Cauchy. 
Since $C(K)^*$ does not contain a copy of $c_0$, 
it is actually unconditionally norm convergent by the 
Bessaga-Pe{\l}czy\'nski Theorem. 
This implies that for the bounded sequence of functions $%
(f_i)_{i\in \N}$ satisfying $f_{i}\to \chi_{G}$ pointwise
constructed in the proof of Lemma~\ref{aux1}, we have 
\begin{eqnarray*}
T^{*}w^{*}(G) &=& \lim_{i\rightarrow \infty } T^{*}w^{*}(f_i) =
\lim_{i\rightarrow \infty } 
\sum_{n=1}^\infty T_n^{*}w^{*}(f_i) \allowdisplaybreaks \\
&=&
\sum_{n=1}^\infty T_n^{*}w^{*}(\chi_{G}) =
\sum_{n=1}^\infty T_n^{*}w^{*}(G)=0.
\end{eqnarray*}
The proof is complete. 
\qed



\end{document}